\documentclass[a4paper]{amsart}

\newtheorem{theorem}{Theorem}

\newtheorem{lemma}{Lemma}

\theoremstyle{definition}

\newtheorem{remark}{Remark}

\usepackage{graphicx}
\usepackage{amssymb,statmath}
\usepackage{mathrsfs}
\newcommand{\R}{{\mathbb{R}}}
\renewcommand{\H}{{\mathbb{H}}}
\newcommand{\<}{\langle}
\renewcommand{\>}{\rangle}
\newcommand{\xii}{{\mathrm{|\xi|}}}
\newcommand{\e}{\varepsilon}
\newcommand{\lin}{\mathrm{lin}}
\newcommand{\lloc}{\mathrm{loc}}
\renewcommand{\L}{{\mathcal{L}}}
\newcommand{\Q}{{\mathcal{Q}}}
\begin{document}

\title[Critical semilinear damped wave equations]{Semilinear damped wave equations with data from Sobolev spaces of negative order: the critical case in Euclidean setting and in the Heisenberg space
}


\author{Marcello D'Abbicco}


\address{M. D'Abbicco, Department of Mathematics, University of Bari, Via E. Orabona 4, 70125 BARI - ITALY}

\email{marcello.dabbicco@uniba.it} 



\maketitle

\begin{abstract}
In this note, we prove the global existence of solutions to the semilinear damped wave equation in $\mathbb{R}^n$, $n\leq6$, with critical nonlinearity under the assumption that the initial data are small in the energy space $H^1\times L^2$ and under the vanishing condition that the initial data belong to $\dot H^{-\gamma}$ for some $\gamma\in(0,n/2)$. A similar result also applies to the damped wave equation in the Heisenberg group $\mathbb{H}^n$, with $n=1,2$.
\end{abstract}

\section{Introduction}
\label{intro}

In this note, we prove that global-in-time solutions to the damped wave equation
\begin{equation}\label{eq:CP} \begin{cases}
u_{tt}-\Delta u +u_t = f(u),& x\in\R^n,\ t>0,\\
(u,u_t)(0,x)=(u_0,u_1), & x\in\R^n,
\end{cases}
\end{equation}
with small initial data in the energy space $H^1\times L^2$ and in $\dot H^{-\gamma}$ exist for $f(u)=|u|^p$ (or, more in general, for $f$ as in~\eqref{eq:f}) in the critical case
\begin{equation}\label{eq:pcrit}
p=1+\frac4{n+2\gamma}, \qquad \gamma\in(0,n/2),
\end{equation}
in space dimension $n=1,2$ and, with some restrictions on $\gamma$, in space dimension $n=3,4,5,6$. This answers the question posed in the recent paper of W. Chen and M. Reissig~\cite[Remark 5]{CR23}, who proved the global existence of small data solutions in the supercritical case, $f(u)=|u|^q$ with $q>p$ (see~\cite[Theorem 1]{CR23}). Here
\[ H^s\cap \dot H^{-\gamma} = \{ u\in H^s: \xii^{-\gamma}\hat u \in L^2 \}, \]
equipped with norm
\[ \|u\|_{H^s\cap \dot H^{-\gamma}} = \|(\<\xi\>^s+\xii^{-\gamma})\hat u\|_{L^2}, \]
where $\hat u=\mathscr{F}u(t,\cdot)$ is the Fourier transform with respect to $x$. The assumption of initial data in homogeneous potential spaces of negative order $\dot H^{-\gamma,p}$ is sometimes called \emph{vanishing condition}; it plays a crucial role in several models, as in the Boussinesq equation~\cite{CO07}.
\begin{theorem}\label{thm:main}
Let $n\leq6$ and let assume
\begin{itemize}
    \item that $\gamma\in(0,n/2)$ if $n=1,2$;
    \item that $\gamma\in(0,\tilde\gamma)$ if $n=3,4$;
    \item that $\gamma\in((n/2)-2,\tilde\gamma)$ if $n=5,6$,
\end{itemize}
where $\tilde\gamma$ is the positive solution to
\[ 2\gamma^2 + n\gamma-2n = 0, \]
if $n\geq3$. Assume that
\begin{equation}\label{eq:f}
f(0)=0, \quad |f(u)-f(v)|\leq C\,|u-v|\big(|u|^{p-1}+|v|^{p-1}\big),
\end{equation}
and $p$ is as in~\eqref{eq:pcrit}. Then there exists $\e_0>0$ such that for any
\begin{equation}\label{eq:dataH1}\begin{split}
u_0\in H^1\cap \dot H^{-\gamma},\quad u_1\in L^2\cap \dot H^{-\gamma},\quad \text{with}\\
\e=\|u_0\|_{H^1\cap\dot H^{-\gamma}}+\|u_1\|_{L^2\cap\dot H^{-\gamma}}\leq\e_0,
\end{split}\end{equation}
there exists a unique solution $u\in\mathcal C([0,\infty),H^1)$ to~\eqref{eq:CP}, where $f$ verifies
 Moreover,
\begin{equation}\label{eq:decay}\begin{split}
\|\nabla u(t,\cdot)\|_{L^2}& \leq C\,\e\,(1+t)^{-\frac{\gamma+1}2},\\
\|u(t,\cdot)\|_{L^2}& \leq C\,\e\,(1+t)^{-\frac\gamma2},
\end{split}\end{equation}
with $C$ independent on $t$ and $\e$.
\end{theorem}
In Theorem~\ref{thm:main} and in the following, $C,C',C''$ denote generic constants which values change in different equations.
\begin{remark}
By dual Sobolev's embeddings (or Hardy-Littlewood-Sobolev inequality), $L^m\cap H^s\hookrightarrow H^s\cap \dot H^{-\gamma}$, for some $m\in(1,2)$, where
\begin{equation}\label{eq:Lm}
\frac1{m} = \frac12+\frac\gamma{n}.
\end{equation}
Therefore, the critical exponent in~\eqref{eq:pcrit} is a generalization of the second critical exponent, in the sense of Lee and Ni~\cite{LN}, derived assuming $L^m$ initial data with $m>1$. Indeed,
\[ 1+\frac4{n+2\gamma}=
1+\frac{2m}n. \]
The global existence of small data solutions for supercritical powers $p>1+2m/n$ for the damped wave equation~\eqref{eq:CP} has been proved in~\cite{IO02}, while the global existence of small data solutions in the critical case~$p=1+2m/n$ has been proved in~\cite{IIOW}. The nonexistence of global solutions for subcritical powers $p<1+2m/n$ has been proved in~\cite[Theorem 3]{DAE17NA}.

We mention that when $m=1$, the critical exponent is Fujita exponent $1+2/n$, and nonexistence of global solutions holds in the critical case, on the contrary of what happens when $m\in(1,2]$, under a suitable sign assumption on the data~\cite{Z01}. On the other hand, there are several other models in which the critical exponent belongs to the existence range, even when $L^1$ data are considered, see, for instance, \cite{DAEP2019}. A sharp integral condition on $f(u)$ in~\eqref{eq:CP} which determines whether global existence of small data solutions holds, or nonexistence of global solutions follows under suitable data sign assumption, has been obtained for the $L^1$ theory in~\cite{EGR2020}.
\end{remark}
The counterpart of nonexistence of global solution for subcritical powers is given by Theorem~3 in~\cite{DAE17NA}, for the more general model of damped $\sigma$-evolution equations; in~\cite{DAE17NA} it is proved that even weak solution $u\in L^p_\lloc$ to~\eqref{eq:CP} with $f(u)=|u|^q$, do not exist when $q\in(1,p)$, assuming a growth condition on the data
\[ u_0\equiv0, \qquad u_1(x) \geq \e\,(1+|x|)^{-\mu}, \quad 1<q<1+\frac2{\mu-1}\,; \]
for any given $q<p$, it is sufficient to choose
\[ 1+\frac2q <\mu < \gamma + \frac{n}2, \]
to obtain the nonexistence of global solutions for any $q<p$. We stress that the growth condition with $\mu$ as above is compatible with the assumption $u\in L^m\hookrightarrow \dot H^{-\gamma}$, namely, $m\mu>n$, where $m$ is as in~\eqref{eq:Lm}. The lifespan estimate
\[ T\leq C\,\e^{-\frac2{2(p'-1)-\mu}} \]
is also provided in~\cite[Theorem 3]{DAE17NA}. A variant of this result with an additional logarithmic term in the growth condition is proved in~\cite{CR23}; also the sharpness of the lifespan estimate is proved in~\cite[Theorem 3]{CR23}.
\begin{remark}
We stress that $H^1\hookrightarrow L^{2p}$ thanks to the assumption %
%
%
$2\gamma\geq n-4$, when $n=5,6$. For $n\geq3$, we may compute
\[ \tilde\gamma = \frac{\sqrt{n^2+16n}-n}4. \]
and notice that $(n/2)-2<\tilde\gamma<(n/2)$ if $n=3,4,5,6$.

These restrictions are consistent with the analogous restrictions related to the use of initial data in $L^m$, where $m$ is as in~\eqref{eq:Lm}, see~\cite{IO02}.
\end{remark}
Theorem~\ref{thm:main} and, similarly, Theorem~1 in~\cite{CR23} may be likely extended to the case of damped $\sigma$-evolution equations, studied in~\cite{DAE17NA}.

Moreover, with a strategy similar to the one employed in the Euclidean case, we may study the damped wave equation on the Heisenberg group~$\H^n$. Namely, we may prove the global existence of solutions to the damped wave equation
\begin{equation}\label{eq:CPH} \begin{cases}
u_{tt}-\Delta_H u +u_t = f(u), & g\in\H^n, \ t>0,\\
(u,u_t)(0,g)=(u_0,u_1)(g), & g\in\H^n,
\end{cases}
\end{equation}
assuming small initial data in $H_\L^1(\H^n)\times L^2(\H^n)$ and in $\dot H_\L^{-\gamma}(\H^n)$, in the critical case
\begin{equation}\label{eq:pcritH}
p=1+\frac4{\Q+2\gamma}, \qquad \gamma\in(0,\tilde\gamma),
\end{equation}
where $\Q=2n+2$ represents the homogeneous dimension of $\H^n$, $n=1,2$, and~$\tilde\gamma$ is the positive solution to
\[ 2\gamma^2 + \Q\gamma-2\Q = 0. \]
Here $\Delta_H=\L$ and $\nabla_H$ denote the sub-Laplacian and the horizontal gradient (we use notation in~\cite{P20}) on the Heisenberg group~$\H^n$, $H_\L^1(\H^n)$ is the Sobolev space of $L^2(\H^n)$ functions~$u$ such that $\nabla_H u\in L^2(\H^n)$, and $H_\L^1(\H^n)\cap \dot H_\L^{-\gamma}(\H^n)$ is the space of $H_\L^1(\H^n)$ functions~$u$ such that the distribution $(-\L)^{\frac{\gamma}2} u$ belongs to $L^2(\H^n)$. Similarly for $L^2(\H^n)\cap \dot H_\L^{-\gamma}(\H^n)$. We use the notation in~\cite{DKMR24} for $(-\L)^{\frac{\gamma}2}$ and for the Sobolev spaces $H_\L^1(\H^n)\cap \dot H_\L^{-\gamma}(\H^n)$ and $L^2(\H^n)\cap \dot H_\L^{-\gamma}(\H^n)$.

This answers the question in~\cite{DKMR24} whether global existence holds or not in the critical case; indeed, in~\cite{DKMR24} both the global existence of small data solutions in the supercritical case, and the nonexistence of global solutions in the subcritical case (under a growth assumption analogous to the one in~\cite{CR23,DAE17NA}), are proved. We address the reader to~\cite{DKMR24,GP20,P20} for more insights about the damped wave equation on the Heisenberg group~$\H^n$.
\begin{theorem}\label{thm:mainH}
Let $n=1,2$ and let assume
\begin{itemize}
    \item that $\gamma\in(0,\tilde\gamma)$ if $n=1$;
    \item that $\gamma\in[1,\tilde\gamma)$ if $n=2$.
\end{itemize}
Then there exists $\e_0>0$ such that for any
\begin{equation}\label{eq:dataHeis}\begin{split}
u_0\in H_\L^1(\H^n)\cap \dot H_\L^{-\gamma}(\H^n),\quad u_1\in L^2(\H^n)\cap \dot H_\L^{-\gamma}(\H^n),\quad \text{with}\\
\e=\|u_0\|_{H_\L^1(\H^n)\cap\dot H_\L^{-\gamma}(\H^n)}+\|u_1\|_{L^2(\H^n)\cap\dot H_\L^{-\gamma}(\H^n)}\leq\e_0,
\end{split}\end{equation}
there exists a unique solution $u\in\mathcal C([0,\infty),H_\L^1(\H^n))$ to~\eqref{eq:CP}, where $f$ verifies~\eqref{eq:f} and $p$ is as in~\eqref{eq:pcritH}. Moreover,
\begin{equation}\label{eq:decayH}\begin{split}
\|\nabla_\L u(t,\cdot)\|_{L^2(\H^n)}& \leq C\,(1+t)^{-\frac{\gamma+1}2},\\
\|u(t,\cdot)\|_{L^2(\H^n)}& \leq C\,(1+t)^{-\frac\gamma2}.
\end{split}\end{equation}
\end{theorem}
The proof of Theorem~\ref{thm:mainH} is very similar to the proof of Theorem~\ref{thm:main}; in particular, we make use of the following result, which is an immediate consequence of Theorem~1.1 in~\cite{DKMR24}, and extends the estimates obtained for the solution to~\eqref{eq:CPH} with $f=0$ in~\cite[Theorem 1.1]{P20}, to the case of initial data in $H_\L^1(\H^n)\cap L^m(\H^n)$.
\begin{lemma}\label{lem:H}
Assume that $u_0\in H^1_\L(\H^n)\cap L^m(\H^n)$ and that $u_1\in L^2(\H^n)\cap L^m(\H^n)$ for some $m\in(1,2]$. Then the solution $u\in\mathcal C^1([0,\infty),H^1(\H^n))$ verifies the following estimates
\[\begin{split}
& \|u(t,\cdot)\|_{L^2(\H^n)}\\
    & \qquad \leq C\,(1+t)^{-\Q\left(\frac1m-\frac12\right)}\,\big(\|u_0\|_{L^2(\H^n)\cap L^m(\H^n)}+\|u_1\|_{H_\L^{-1}(\H^n)\cap L^m(\H^n)}\big),\\
& \|\nabla_H u(t,\cdot)\|_{L^2(\H^n)}\\
    & \qquad \leq C\,(1+t)^{-\Q\left(\frac1m-\frac12\right)-\frac12}\,\big(\|u_0\|_{H_\L^1(\H^n)\cap L^m(\H^n)}+\|u_1\|_{L^2(\H^n)\cap L^m(\H^n)}\big).
\end{split}\]
\end{lemma}
In a similar way to the proof of Theorem~\ref{thm:mainH}, we expect that it is possible to prove the global existence of small data solutions in the critical case
\[ p=1+\frac{2\nu}{\Q+2\gamma} \]
for higher order hypoelliptic damped wave equations on graded lie groups $\mathbb{G}$, with small data from negative order Sobolev spaces recently studied in~\cite{DKMR24b},
\[ \begin{cases}
u_{tt}+\mathcal R u +u_t = f(u), & x\in\mathbb{G}, \ t>0,\\
(u,u_t)(0,x)=(u_0,u_1)(x), & x\in\mathbb{G},
\end{cases} \]
where $\mathcal R$ is a Rockland operator of homogeneous degree~$\nu\geq2$.

\section{Proof of Theorem~\ref{thm:main}}

We now prove Theorem~\ref{thm:main}.
\begin{proof}
For any $T>0$, let $X$ be the Banach space of functions in $\mathcal C([0,T],H^1)$ equipped with the norm
\[ \|u\|_X=\sup_{t\in[0,T]} t^{\frac\gamma2} \big( \|u(t,\cdot)\|_{L^2} + t^{\frac12}\,\|\nabla u(t,\cdot)\|_{L^2}\big). \]
Thanks to Proposition~1 in~\cite{CR23}, the solution
\[ u_\lin=K_0(t)\ast u_0 + K_1(t)\ast u_1\]
to the linear problem~\eqref{eq:CP} with $f=0$ is in $X$ and
\begin{equation}\label{eq:ueps}
\|u\|_X\leq C_1\varepsilon,
\end{equation}
for some $C_1>0$, independent on $T$, where $\e$ is as in~\eqref{eq:dataH1}. Thanks to Gagliardo-Nirenberg inequality, for any $u\in X$ it holds
\begin{equation}\label{eq:uq}
\|u(t,\cdot)\|_{L^q} \leq (1+t)^{-\frac{n}2\left(\frac12-\frac1q\right)-\frac\gamma2} \|u\|_X,
\end{equation}
for any $q\in[2,\infty]$ if $n=1$, for any $q\in[2,\infty)$ if $n=2$ and for any $q\in[2,2^*]$, where
\[ \frac1{2^*}=\frac12-\frac1n, \quad \text{i.e.}\qquad 2^*=\frac{2n}{n-2}, \]
if $n\geq3$. We now consider the Duhamel's integral operator
\[ Nu(t,\cdot) = \int_0^t K_1(t-\tau)\ast f(u(\tau,\cdot))\,d\tau. \]
We claim that for any $u,v\in X$, it holds
\begin{equation}\label{eq:contr}
\|Nu-Nv\|_X \leq C_2 D(u,v),\qquad D(u,v)=\|u-v\|_X \big(\|u\|_X^{p-1}+\|v\|_X^{p-1}\big),
\end{equation}
for some~$C_2>0$, independent on $T$. To prove~\eqref{eq:contr}, we use the well-known estimates (see Lemma 1 in~\cite{M76}):
\[ \|\nabla^j K_1(t)\ast g\|_{L^q} \leq C\,(1+t)^{-n\left(\frac1m-\frac1q\right)-\frac{j}2}\,\|g\|_{L^m\cap L^2},\quad \|g\|_{L^m\cap L^2}=\|g\|_{L^m}+\|g\|_{L^2}, \]
that hold for $m\in[1,2]$ and $j=0,1$. For $j=0,1$, we get
\[ \begin{split}
& \|\nabla^j (Nu-Nv)(t,\cdot)\|_{L^2}\\
     & \qquad \leq \int_0^{t/2} (1+t-\tau)^{-\frac{n}2\left(\frac1m-\frac12\right)-\frac{j}2} \, \|(f(u)-f(v))(\tau,\cdot)\|_{L^m\cap L^2}d\tau \\
     & \qquad +\int_{t/2}^t (1+t-\tau)^{-\frac{j}2} \|(f(u)-f(v))(\tau,\cdot)\|_{L^2}\,d\tau,
\end{split} \]
where $m=1$ if $p\geq2$ or $m=2/p$ otherwise (we use $m=2$ in $[t/2,t]$). Due to~\eqref{eq:f} and~\eqref{eq:uq}, we obtain
\[ \begin{split}
& \|(f(u)-f(v))(\tau,\cdot)\|_{L^2}\\
    & \qquad \leq C\,\|(u-v)(\tau,\cdot)\|_{L^{2p}}\big(\|u(\tau,\cdot)\|_{L^{2p}}^{p-1}+\|v(\tau,\cdot)\|_{L^{2p}}^{p-1}\big) \\
    & \qquad\leq C'(1+\tau)^{-\frac{n}4(p-1)-\frac\gamma2p}\,D(u,v) = C'(1+\tau)^{-1-\frac\gamma2}\,D(u,v),
    \end{split} \]
so that we find
\[ \begin{split}
& \int_{t/2}^t (1+t-\tau)^{-\frac{j}2} \|(f(u)-f(v))(\tau,\cdot)\|_{L^2}\,d\tau \\
& \qquad \leq C \int_{t/2}^t (1+t-\tau)^{-\frac{j}2} (1+\tau)^{-1-\frac\gamma2}\,d\tau\,D(u,v)\\
& \qquad \leq C' (1+t)^{-1-\frac\gamma2}\int_{t/2}^t (1+t-\tau)^{-\frac{j}2} \,d\tau\,D(u,v)\\
& \qquad \leq C''\,(1+t)^{-\frac{\gamma+j}2}\,D(u,v).
\end{split} \]
We now consider the first integral. First, let $p\geq2$, so that $m=1$. In this case, due to~\eqref{eq:f} and~\eqref{eq:uq}, we obtain
\[ \begin{split}
& \|(f(u)-f(v))(\tau,\cdot)\|_{L^1}\\
    & \qquad \leq C\,\,\|(u-v)(\tau,\cdot)\|_{L^p}\big(\|u(\tau,\cdot)\|_{L^p}^{p-1}+\|v(\tau,\cdot)\|_{L^p}^{p-1}\big)\\
    & \qquad \leq C'(1+\tau)^{-\frac{n}4(p-2)-\frac\gamma2p} D(u,v) = C'(1+\tau)^{\frac{n}4-1-\frac\gamma2} D(u,v),
\end{split} \]
so that we get
\[ \begin{split}
& \int_0^{t/2} (1+t-\tau)^{-\frac{n}4-\frac{j}2} \, \big(\|(f(u)-f(v))(\tau,\cdot)\|_{L^1\cap L^2}\big)d\tau \\
& \qquad \leq C\,(1+t)^{-\frac{n}4-\frac{j}2} \,\int_0^{t/2} (1+\tau)^{\frac{n}4-1-\frac\gamma2}\, d\tau\,D(u,v)\\
& \qquad \leq C'\,(1+t)^{-\frac{\gamma+j}2}\,D(u,v).
\end{split} \]
where we used that
\[ \frac{n}4-1-\frac\gamma2>-1, \]
as a consequence of $\gamma<n/2$, to estimate the integral. Now let $p\leq2$ so that $m=2/p$. In this case, due to~\eqref{eq:f} and~\eqref{eq:uq}, we obtain
\[ \begin{split}
& \|(f(u)-f(v))(\tau,\cdot)\|_{L^{\frac2p}}\\
    & \qquad \leq C\,\|(u-v)(\tau,\cdot)\|_{L^2}\big(\|u(\tau,\cdot)\|_{L^2}^{p-1}+\|v(\tau,\cdot)\|_{L^2}^{p-1}\big) \\
    & \qquad\leq C'(1+\tau)^{-\frac\gamma2p}\,D(u,v),
\end{split} \]
so that we get
\[ \begin{split}
& \int_0^{t/2} (1+t-\tau)^{-\frac{n}4(p-1)-\frac{j}2} \, \|f(u)(\tau,\cdot)\|_{L^{\frac2m}\cap L^2}d\tau \\
& \qquad \leq C\,(1+t)^{-\frac{n}4(p-1)-\frac{j}2} \,\int_0^{t/2} (1+\tau)^{-\frac\gamma2p}\, d\tau\,\|u\|_X^p\\
& \qquad \leq C'\,(1+t)^{-\frac{\gamma+j}2}\,\|u\|_X^p,
\end{split} \]
where we used that
\[ \frac\gamma2p = \frac\gamma2 \left(1+\frac4{n+2\gamma}\right)<1, \]
as a consequence of the assumption %
%
%
$\gamma<\tilde\gamma$ when $n\geq3$, to estimate the integral.

This concludes the proof of~\eqref{eq:contr}.

Let
\begin{equation}\label{eq:R}
R = 2C_1\,\e,
\end{equation}
where $C_1$ is as in~\eqref{eq:ueps}, and fix~$\e_0$ such that $C_2 R^{p-1}\leq 1/4$, for any~$\e\leq\e_0$, where $C_2$ is as in~\eqref{eq:contr}. Then the operator
\[ \Lambda: u\in X \mapsto \Lambda u = u_\lin + Nu \in X \]
maps the closed ball of radius~$R$ in~$X$ into itself, since
\[ \|Nu+u_\lin\|_X \leq \frac14\,\|u\|_X + \|u_\lin\|_X \leq R. \]
Moreover, the operator $\Lambda$ is a contraction, due to
\[ \|\Lambda u-\Lambda v\|_X = \|Nu-Nv\|_X \leq \frac12\,\|u-v\|_X. \]
Therefore, in the closed ball of radius~$R$ in~$X$, there exists a unique fixed point~$u$ for the operator~$\Lambda$, for sufficiently small data, i.e., $\e\leq\e_0$ as in~\eqref{eq:dataH1}. This fixed point is the unique solution to~\eqref{eq:CP} in~$\mathcal C([0,T],H^1)$. Moreover, the fixed point verifies the estimate
\begin{equation}\label{eq:decayX}
\|u\|_X \leq R = 2C_1\,\e.
\end{equation}
Due to the fact that no constant depends on~$T$, this argument gives the existence of a unique global-in-time solution in~$\mathcal C([0,\infty),H^1)$, and~\eqref{eq:decayX} provides the desired decay estimate~\eqref{eq:decay}. This concludes the proof.
\end{proof}
For simplicity, we considered initial data in the energy space $H^1\times L^2$ and in $\dot H^{-\gamma}$, but the statement and the proof may be easily adapted for initial data in $H^s\times H^{s-1}$ and in $\dot H^{-\gamma}$, $s\in(0,1)$.

\section{Proof of Lemma~\ref{lem:H} and of Theorem~\ref{thm:mainH}}

The proof of Theorem~\ref{thm:mainH} is very similar to the proof of Theorem~\ref{thm:main}, so we focus on the main differences. In order to follow the same outline of the proof of Theorem~\ref{thm:main}, we use Lemma~\ref{lem:H}, whose proof follows.
\begin{proof}[Proof of Lemma~\ref{lem:H}]
We fix~$\eta\in(0,\Q/2)$ such that
\[ \frac1m = \frac12 + \frac\eta\Q, \quad \text{i.e.}\qquad \eta=\Q\,\left(\frac1m-\frac12\right); \]
hence, by Hardy–Littlewood–Sobolev inequality on the Heisenberg group (see, for instance, Theorem 2.1 in~\cite{DKMR24}), we get
\begin{equation}\label{eq:dualH}
\|u(t,\cdot)\|_{L^2\cap \dot H_\L^{-\eta}} \leq C\,\|u(t,\cdot)\|_{L^m\cap L^2}.
\end{equation}
By Theorem~1.1 in~\cite{DKMR24}, thanks to~\eqref{eq:dualH}, we obtain
\[ \begin{split}
& \|\nabla_H u(t,\cdot)\|_{L^2}\leq C\,\|u(t,\cdot)\|_{\dot H^1_\L} \\
    & \qquad \leq C'\,(1+t)^{-\frac{\eta+1}2}\,\big(\|u_0\|_{H_\L^1(\H^n)\cap \dot H_\L^{-\eta}}+\|u_1\|_{L^2(\H^n)\cap \dot H_\L^{-\eta}}\big)\\
    & \qquad \leq C''\,(1+t)^{-\frac{\Q}2\left(\frac1m-\frac12\right)-\frac12}\,\big(\|u_0\|_{H_\L^1(\H^n)\cap L^m(\H^n)}+\|u_1\|_{L^2(\H^n)\cap L^m(\H^n)}\big),\\
& \|u(t,\cdot)\|_{L^2} \\
    & \qquad \leq C\,(1+t)^{-\frac\eta2}\,\big(\|u_0\|_{L^2(\H^n)\cap \dot H_\L^{-\eta}}+\|u_1\|_{H^{-1}_\L(\H^n)\cap \dot H_\L^{-\eta}}\big)\\
    & \qquad \leq C'\,(1+t)^{-\frac{\Q}2\left(\frac1m-\frac12\right)}\,\big(\|u_0\|_{L^2(\H^n)\cap L^m(\H^n)}+\|u_1\|_{H^{-1}_\L(\H^n)\cap L^m(\H^n)}\big)
\end{split} \]
This concludes the proof.
\end{proof}
We now prove Theorem~\ref{thm:mainH}.
\begin{proof}[Proof of Theorem~\ref{thm:mainH}]
Thanks to Theorem~1.1 in~\cite{DKMR24}, the solution $u_\lin$ to~\eqref{eq:CPH} with $f=0$ belongs to the solution space $X=\mathcal C^1([0,T],H_\L^1)$ and~\eqref{eq:ueps} holds for some $C_1>0$, independent on $T$, where $\e$ is as in~\eqref{eq:dataHeis}, where
\[ \|u\|_X=\sup_{t\in[0,T]} t^{\frac\gamma2} \big( \|u(t,\cdot)\|_{L^2} + t^{\frac12}\,\|\nabla_H u(t,\cdot)\|_{L^2}\big). \]
For brevity, in the formula above and in the following we write $L^2$ in place of $L^2(\H^n)$ and similarly for the Sobolev spaces.

Thanks to Gagliardo-Nirenberg inequality (see~\cite{CR13}, see also~\cite[Lemma 4.1]{GP20}), for any $u\in X$ it holds
\begin{equation}\label{eq:uqH}
\|u(t,\cdot)\|_{L^q} \leq (1+t)^{-\frac{\Q}2\left(\frac12-\frac1q\right)-\frac\gamma2} \|u\|_X,\qquad q\in[2,2^*],
\end{equation}
where
\[ \frac1{2^*}=\frac12-\frac1\Q, \quad \text{i.e.}\qquad 2^*=\frac{2\Q}{\Q-2}=2+\frac2n. \]
We claim that for any $u,v\in X$, \eqref{eq:contr} holds for some~$C_2>0$, independent on $T$, where $N$ is the Duhamel's integral operator associated to~\eqref{eq:CPH}, i.e.,
\[ Nu(t,\cdot) = \int_0^t E_1(t-\tau)\ast_{(g)} f(u(\tau,\cdot))\,ds,\]
where now $E_1$ is the fundamental solution to~\eqref{eq:CPH} and $\ast_{(g)}$ is the group convolution product on $\H^n$ with respect to the $g$ variable (see~\cite{GP20}).

To prove~\eqref{eq:contr}, we proceed as in the proof of Theorem~\ref{thm:main}. Noticing that $p\in(1,2)$, so that $m=2/p\in(1,2)$, by Lemma~\ref{lem:H}, for $j=0,1$, we obtain
\[ \begin{split}
& \|\nabla_H^j (Nu-Nv)(t,\cdot)\|_{L^2}\\
     & \qquad \leq \int_0^{t/2} (1+t-\tau)^{-\frac{\Q}4(p-1)-\frac{j}2} \, \|(f(u)-f(v))(\tau,\cdot)\|_{L^2\cap L^{\frac2p}}d\tau \\
     & \qquad +\int_{t/2}^t (1+t-\tau)^{-\frac{j}2} \|(f(u)-f(v))(\tau,\cdot)\|_{L^2}\,d\tau,
\end{split} \]
as in the proof of Theorem~\ref{thm:main}, with $\Q$ in place of $n$. Due to~\eqref{eq:f} and~\eqref{eq:uqH}, we obtain
\[ \begin{split}
& \|(f(u)-f(v))(\tau,\cdot)\|_{L^2}\\
    & \qquad \leq C\,\|(u-v)(\tau,\cdot)\|_{L^{2p}}\big(\|u(\tau,\cdot)\|_{L^{2p}}^{p-1}+\|v(\tau,\cdot)\|_{L^{2p}}^{p-1}\big) \\
    & \qquad\leq C(1+\tau)^{-\frac{\Q}4(p-1)-\frac\gamma2p}\,D(u,v) = C(1+\tau)^{-1-\frac\gamma2}\,D(u,v),
    \end{split} \]
so that we find
\[ 
\int_{t/2}^t (1+t-\tau)^{-\frac{j}2} \|(f(u)-f(v))(\tau,\cdot)\|_{L^2}\,d\tau 
\leq C\,(1+t)^{-\frac{\gamma+j}2}\,D(u,v),
\]
as in the proof of Theorem~\ref{thm:main}. We now consider the first integral. Thanks to~\eqref{eq:dualH}, we may proceed as in the proof of Theorem~\ref{thm:main}; due to~\eqref{eq:f} and~\eqref{eq:uqH}, we obtain
\[ \begin{split}
& \|(f(u)-f(v))(\tau,\cdot)\|_{L^{\frac2p}} \\
    & \qquad \leq C\,\|(u-v)(\tau,\cdot)\|_{L^2}\big(\|u(\tau,\cdot)\|_{L^2}^{p-1}+\|v(\tau,\cdot)\|_{L^2}^{p-1}\big) \\
    & \qquad\leq C'(1+\tau)^{-\frac\gamma2p}\,D(u,v),
\end{split} \]
so that we get
\[ 
\int_0^{t/2} (1+t-\tau)^{-\frac{\Q}4(p-1)-\frac{j}2} \, \|(f(u)-f(v))(\tau,\cdot)\|_{L^{\frac2p}\cap L^2}d\tau 
\leq C\,(1+t)^{-\frac{\gamma+j}2}\,D(u,v),
\]
where we used $\gamma<\tilde\gamma$ as in the proof of Theorem~\ref{thm:main}. This concludes the proof of~\eqref{eq:contr}.

The proof follows as in the last part of the proof of Theorem~\ref{thm:main}.
\end{proof}

%

\section*{Acknowledgement}

The author is supported by PNR MUR Project CUP\_H93C22000450007, and by INdAM-GNAMPA Project CUP\_E55F22000270001. The author thanks Prof. Palmieri for useful discussions on the Heisenberg group.



\end{document}